\documentclass[12pt]{amsart}

\usepackage{amsmath,amsthm,amsfonts,latexsym,amscd,amssymb,enumerate}
\usepackage{amscd}
\usepackage[hypertex,backref]{hyperref}

\newtheorem{theorem}{Theorem}[section]
\newtheorem{thm}[theorem]{Theorem}
\newtheorem{prop}[theorem]{Proposition}
\newtheorem{lem}[theorem]{Lemma}

\theoremstyle{remark}
 
\newtheorem{rem}[theorem]{Remark}

\theoremstyle{definition}
\newtheorem{defn}[theorem]{Definition}

\theoremstyle{definition}

\theoremstyle{remark}

\DeclareMathOperator{\ind}{ind}

\DeclareMathOperator{\End}{End}
\DeclareMathOperator{\ch}{{\rm ch}}
\DeclareMathOperator{\tr}{tr}

\DeclareMathOperator{\Mat}{Mat}

\newcommand{\K}{\mathbb{ K}}
\newcommand{\C}{\mathbb{ C}}

\newcommand{\N}{\mathbb{ N}}

\newcommand{\Z}{\mathbb{ Z}}

\newcommand{\Hom}{\operatorname{Hom}}

\newcommand{\Q}{\mathbb{ Q}}
\newcommand{\R}{\mathbb{ R}}

{\catcode`@=11\global\let\c@equation=\c@theorem}

\begin{document}

\title{Enlargeability and index theory: Infinite covers}

\author{Bernhard Hanke and Thomas Schick}

\begin{abstract} In \cite{enlind} we showed non-vanishing of the universal index elements in 
the $K$-theory of the maximal $C^*$-algebras of the fundamental groups of enlargeable spin manifolds. 
The underlying notion of enlargeability was the one from \cite{GL}, 
involving contracting maps defined on finite covers of the given manifolds. In the paper at hand, 
we weaken this assumption to the one in \cite{GLIHES} where
infinite covers are allowed. The new idea is the construction of a
geometrically given $C^*$-algebra with trace which encodes the information
given by these infinite covers. Along the way we obtain an easy proof of a
relative index theorem relevant in this context.
\end{abstract}

\address{Georg-August-Universit{\"a}t G\"ottingen \\ Germany \\}
\email{hanke@uni-math.gwdg.de \vspace{1cm}}

\address{Georg-August-Universit{\"a}t G{\"o}ttingen  \\
Germany \\}
\email{schick@uni-math.gwdg.de}

\thanks{We thank S. Stolz and A. Thom for useful conversations regarding the 
research in this paper. Both authors are members of the 
DFG emphasis programme ``Globale Differentialgeometrie'' whose 
support is gratefully acknowledged.}

\maketitle

\section{Overview of results}

This paper is a sequel to \cite{enlind}. We quickly review the relevant mathematics and 
the motivation for our research. 
Since the work of Lichnerowicz it has been known that index theory can be used to construct 
computable obstructions to the existence of positive scalar curvature metrics on closed spin manifolds $M^n$: 
If such a metric exists on $M$, then $\hat{A}(M)$, a characteristic number of $M$ involving its 
Pontryagin classes, vanishes. However, the converse is not true: In \cite{GL} it was shown
(again by index theoretic methods) that 
also enlargeability (see the definition below) of $M$ is an obstruction to the existence of a 
positive scalar curvature metric. Furthermore, it is well known that there are enlargeable manifolds 
with vanishing $\hat{A}$-genus (for example tori). The $\hat{A}$-obstruction was 
gradually refined and in its most sophisticated form constructed by Rosenberg \cite{Ros} as an element 
\[    
   \alpha^{\R}_{max}(M) \in KO_n(C^*_{max, \R} \pi_1(M)) 
\]
in the real $K$-theory of the maximal $C^*$-algebra of $\pi_1(M)$. It
is often assumed that this obstruction 
subsumes all other index theoretic obstructions against positive scalar
curvature metrics on $M$. 
Results of Stolz \cite{Stolz} imply that this is indeed the case, if the Baum-Connes 
conjecture holds for $\pi_1(M)$. Note that the emphasis here is on index theoretic
obstructions. The example of \cite{Schick} shows that in general there are
other obstructions which can not be detected by $\alpha^{\R}_{max}(M)$.

Guided by these facts, it is natural to ask whether non-vanishing of $\alpha^{\R}_{max}(M)$ for enlargeable 
$M$ can be shown independently of the Baum-Connes conjecture. In our paper \cite{enlind} this was achieved 
for manifolds that are enlargeable in the sense of \cite{GL}. It is the purpose of this 
note to extend this result to manifolds that are enlargeable in the more general 
sense of \cite{GLIHES}. 

Let us start with the relevant definition introduced in \cite{GLIHES}. 

\begin{defn} \label{general} We call a connected closed oriented manifold $M^n$ {\it enlargeable} if the following 
holds: Fix some Riemannian metric $g$ on $M$. Then, for all $\epsilon > 0$, there 
is a connected covering  $\overline{M} \to M$ (which may be finite or infinite) and an $\epsilon$-contracting 
map $(\overline{M}, \overline{g}) \to (S^n, g_0)$ which is constant outside a compact subset 
of $\overline{M}$ and of nonzero degree. Here, $\overline{g}$ is induced by $g$ and $g_0$ is 
the standard metric on $S^n$. 

The manifold $M$ is called {\em area-enlargeable} if in the above definition $\epsilon$-contracting 
is replaced by $\epsilon$-area-contracting (cf.~\cite{enlind}, Def. 1.1). 
\end{defn}

If we require in addition the covers $\overline{M} \to M$ to be finite, this coincides with 
Definition 1.1 in \cite{enlind} which goes back to \cite{GL}. It is an easy consequence of 
Hadamard's theorem that $M$ is enlargeable, if 
it admits a Riemannian metric of nonnegative sectional curvature. Moreover, the 
connected sum of an enlargeable manifold with an arbitrary connected closed 
oriented manifold is again enlargeable. 

>From now on the notion of (area-)enlargeability refers to Definition \ref{general}. In contrast, manifolds 
sharing the more restrictive property used in \cite{enlind} will be called {\it compactly} 
(area-)enlargeable.   

It is obvious that each enlargeable manifold is area-enlargeable, but  
not known if the converse holds. 

Because we are mainly interested in explaining the method, we confine
ourselves
to indicating a proof of the following analogue of one of the main results 
of \cite{enlind}.

\begin{thm} \label{main1} Let $M$ be an enlargeable or area-enlargeable spin manifold of 
dimension $n$. Then 
\[
   \alpha_{max}(M) \neq 0 \in K_n(C^*_{max} \pi_1(M)) \, . 
\]
\end{thm}

Here and in the rest of this paper we use complex $K$-theory. The statement of Theorem 
\ref{main1} implies a similar non-vanishing statement for $\alpha^{\R}_{max}(M)$. 

We mention the following related result. The proof of the corresponding statement in \cite{enlind} 
applies without change. 

\begin{thm} Let $M$ be an area-enlargeable manifold of dimension $n$ and let 
\[
   f : M \to B \pi_1(M) 
\]
be the classifying map of the universal cover of $M$. Then the image of the fundamental 
class of $M$ under the induced map in homology 
\[
   f_* : H_n(M ; \Q) \to H_n( B \pi_1(M) ; \Q) 
\]
is different from zero. 
\end{thm}

The starting point for the discussion in \cite{enlind} is the observation (inspired by the work 
of Gromov and Lawson) that on an even dimensional compactly area-enlargeable manifold, one can construct an 
{\rm almost flat} vector bundle whose total Chern class is different from 
zero exactly in degrees $0$ and $n$ (where $\dim M = 2n$). 
By definition (cf. \cite{enlind}, Section 2), this means  that we have a sequence $(E_i)_{i \in \N}$ of 
finite dimensional unitary vector bundles on $M$ whose total Chern classes
have the described property and that enjoy the following asymptotic flatness property: 
Each bundle $E_i$ is equipped 
with a metric connection $\nabla_i$ so that the 
associated curvature $2$-forms ($d_i := \dim E_i$)
\[  
       \Omega_i \in \Omega^2(M ; \mathfrak{u}(d_i)) 
\]
vanish asymptotically in the sense that 
\[
    \lim_{i \to \infty} \| \Omega_i \| = 0 \, . 
\]
Here we use the maximum norm on the unit sphere bundle in $\Lambda^2 M$ and the 
operator norm on $\mathfrak{u}(d_i)\subset \Mat(d_i) = \C^{d_i \times d_i}$. 

The main contribution of the paper at hand is the proof of a similar statement 
under the assumption that $M$ is area-enlargeable in the sense of Definition \ref{general}. 
The main difference is that now the constituents of the almost flat bundle
may be infinite dimensional. 

\begin{defn} \label{almflat} Let $M$ be a closed smooth Riemannian manifold. 
An {\em almost flat Hilbert-module bundle} consists of 
a sequence $(C_i)_{i \in \N}$ of complex $C^*$-algebras (which may be unital or non-unital) 
and a sequence $(F_i, \nabla_i)$ of Hilbert $C_i$-module bundles over $M$ whose fibres are $C_i$-isomorphic to 
projective right $C_i$-modules of the form $q_i C_i$, where $q_i \in C_i$ is a projection. 
Furthermore we require that each $F_i$ is equipped with a metric $C_i$-linear 
connection 
\[   
     \nabla_i : \Gamma(F_i) \to \Gamma(T^*M \otimes F_i)   
\]
so that the associated curvature forms
\[
    \Omega_i \in \Omega^2 (M ; \End(F_i)) 
\]
tend to zero in the sense that 
\[ 
    \lim_{i \to \infty} \| \Omega_i \| = 0 \, . 
\]
\end{defn} 

For a detailed discussion of  Hilbert module bundles, we refer to \cite{Sch}. 

The essential step in our discussion is the following technical result. We 
denote by $\K$ the $C^*$-algebra of compact operators on the Hilbert space $l^2(\N)$. 

\begin{prop} \label{technisch} Let $M^{2n}$ be an even dimensional area-enlargeable spin manifold
and let $i \in \N$ be a positive natural number. Then there is a $C^*$-algebra 
$C_i$ and a Hilbert $C_i$-module bundle $F_i \to M$ with connection $\nabla_i$ 
with the following properties: The curvature $\Omega_i$ of $F_i$ satisfies 
\[
    \| \Omega_i \| \leq \frac{1}{i} C 
\]
where $C$ is a constant depending only on $n$. Furthermore, there is a split extension of the form 
\[
     0 \to \K \to C_i \to \Gamma_i \to 0    
\]
with a certain $C^*$-algebra $\Gamma_i$. 
In particular, each $K_0(C_i)$ canonically splits off a $\Z=K_0(\K)$
summand. Let $a_i \in K_0(C_i)$ denote the index of the spin Dirac operator on
$M$ twisted with  
$F_i$. Then the $\Z=K_0(\K)$-component $z_i$ of $a_i$ is different from $0$. 
\end{prop}

We make some remarks concerning the idea behind this statement. Using the area-enlargeability 
of $M$, there is a cover $\overline{M} \to M$ and a unitary bundle $F \to \overline{M}$  
of finite dimension $d$ whose curvature norm is bounded by $\frac{1}{i} C$. Moreover, if $\overline{M}$ is non-compact, 
this bundle is trivial outside a certain compact subset of $\overline M$. The bundle $F$ has the property
that the (relative) 
index of the Dirac operator on $\overline{M}$ twisted with the virtual bundle $F - \underline{\C^d}$
is different from zero. The  
construction of $F$ is explained in \cite{GLIHES} and will be reviewed in Section \ref{Beweis} below. Proposition 
\ref{technisch} claims that this relative index (which will be equal to $z_i$) can be read off 
in a canonical way from the index of a twisted Dirac operator on $M$. This argument is carried out 
in our previous paper \cite{enlind}, if $\overline{M}$ is compact (i.e. if the 
cover $\overline{M} \to M$ is finite). The purpose of Proposition \ref{technisch} is 
a generalization of this construction to the case of non-compact $\overline{M}$.

Before proving Proposition \ref{technisch} in the next section, we explain how Theorem \ref{main1} 
follows. 

By an easy suspension argument (cf. the proof of Theorem 4.2. in \cite{enlind})  
we may restrict attention to even dimensional $M$. 

Let $A$ denote the complex $C^*$-algebra of norm bounded sequences
\[ 
    (a_i)_{i \in \{1,2,\ldots \}} \in \prod_{i=1}^{\infty} C_i 
\]
and for $i \in \N$, we denote by $A_i \subset A$ the subalgebra where all but the $i$th entry vanish. 
This subalgebra can be identified with $C_i$. Let 
\[
    q := (q_i) \in A 
\]
be the projection consisting of the projections $q_i$. We now  
have an analogue of \cite{enlind}, Theorem 2.1: The bundles $F_i$ together with the connections $\nabla_i$ 
can be assembled to a smooth Hilbert $A$-module bundle $V$ together with an $A$-linear metric connection 
\[
     \nabla^V : \Gamma(V) \to \Gamma(T^*M \otimes V) \, 
\]
such that for each $i \in \N$, the subbundle $V_i := V \cdot A_i$ is isomorphic to $F_i$ as 
an $A_i$-Hilbert module bundle, the connection $\nabla^V$ preserves the subbundles $V_i$ and the 
norms of the induced connections $\nabla^V_i$ on $V_i$ tend to zero.   The proof of this statement 
is similar to the proof of \cite[Theorem 2.1]{enlind} and will be omitted.

We now argue as in Section 3 of \cite{enlind} and define 
\[
   A' = \overline{\bigoplus_{i=1}^{\infty}  A_i} \subset A 
\]
as the closed two sided ideal consisting of sequences converging to $0$. Let 
\[
   Q : = A / A'
\]
denote the quotient $C^*$-algebra. Similarly, we define
$A_\Gamma:=\prod_{i\in\N}\Gamma_i$ with ideal
$A'_\Gamma:=\bigoplus_{i\in\N}\Gamma_i$ and quotient
$Q_\Gamma:=A_\Gamma/A'_\Gamma$. The bundle 
\[
   W : = V / (V \cdot A')= V\times_A Q \to M 
\]
is a smooth Hilbert $Q$-module bundle with fibre $\overline{q} Q$ where $\overline{q} \in Q$ 
is the image of $q$. The connection $\nabla^V$ induces a $Q$-linear metric connection 
on $W$. The crucial point is that this connection is flat
(cf. \cite[Proposition 4.3]{enlind}). After fixing 
a basepoint in $x \in M$ and an isomorphism $W_x \cong \overline{q} Q$, we hence get a holonomy representation 
\[ 
    \pi_1(M,x) \to \Hom_Q(\overline{q}Q, \overline{q} Q) = \overline{q}Q \overline{q} 
\]
with values in the unitaries of $\overline{q} Q \overline{q}$. Using the universal property of 
$C^*_{max} \pi_1(M)$, we obtain a $C^*$-algebra homomorphism
\[
     \phi : C^*_{max} \pi_1(M) \to \overline{q} Q \overline{q} \to Q  
\]
and hence an induced map in $K$-theory
\[
   \phi_* =   K_0(C^*_{max} \pi_1(M)) \to K_0(Q)  \, . 
\]
One now shows that the image of $\alpha_{max}(M)$ under this map is different from $0$. For 
that purpose we need some information on the $K$-theory of $Q$. Consider
the following commutative diagram
\begin{equation*}
  \begin{CD}
    && 0 && 0 && 0\\
    && @VVV  @VVV  @VVV\\
    0 @>>> \bigoplus_{i\in\N}\K @>>> A' @>>> A'_\Gamma @>>>
    0\\
    && @VVV  @VVV  @VVV\\
    0 @>>> \prod_{i\in\N} \K @>>> A @>>> A_\Gamma @>>> 0 \\
    && @VVV @VVV @VVV\\
    0 @>>> (\prod_{i\in\N}\K)/(\bigoplus_{i\in\N}\K) @>>> Q @>>>
    Q_\Gamma@>>> 0\\
    && @VVV @VVV @VVV\\
    && 0 && 0 && 0
  \end{CD}
\end{equation*}
where the horizontal rows are split exact (and the splits also make the
diagram commutative) and the vertical columns are exact.
In \cite[Proposition 3.5]{enlind}, we calculated the K-theory of the left
vertical column. It follows that the group $K_0(Q)$ splits off 
a summand 
\[
    \big( \prod_{i=1}^{\infty} \Z \big) / \big( \bigoplus_{i=1}^{\infty} \Z \big) 
\]
which is the image of a corresponding summand $\prod_{i\in\N}\Z$ in
$K_0(A)$ under the map induced by the projection.
Moreover, as in \cite{enlind} one now checks 
that the component of $\phi_*(\alpha_{max}(M))$ in this summand is represented by 
the sequence 
\[ 
    (z_1, z_2, \ldots)  
\] 
and is therefore different from $0$. This implies that $\phi_*(\alpha_{max}(M)) \neq 0$ and hence 
$\alpha_{max}(M) \neq 0$ which completes the proof of Theorem \ref{main1}.

\section{Proof of Proposition \ref{technisch}} \label{Beweis} 

From now on, fix $i\in\N$.
Since $M^{2n}$ is area-enlargeable,  there is a covering 
\[
    f:  \overline{M} \to M 
\]
together with an $\frac{1}{i}$-area-contracting map 
\[
      \psi : \overline{M} \to S^{2n} 
\]
which is constant outside a certain compact subset $K \subset \overline M$ (which will 
be fixed from now on) and of nonzero degree. If $\overline{M}$ is compact, we 
set $K = \overline{M}$. 

We fix a finite dimensional unitary vector bundle 
\[
    E \to S^{2n}
\]
whose $n$-th Chern class is different from zero. This is possible because the Chern character 
is a rational isomorphism. Let $d$ be the complex 
dimension of $E$. We fix a unitary connection on $E$. 

Now let 
\[
      F : = \psi^*(E) \to \overline{M} 
\]
be the pull back bundle with its induced unitary structure and connection
$\nabla_F$.
Let 
\[
    \Omega_F   \in  \Omega^2(\overline{M} ; \mathfrak{u}(d)) \, , ~~ \Omega_E \in \Omega^2(S^n  ; \mathfrak{u}(d))  
\]
be the curvature forms of $F$ and $E$. We set
\[
    C := \| \Omega_E \| \, . 
\]
Because $\psi$ is $\frac{1}{i}$-area 
contracting,  we have 
\[
    \| \Omega_F \| \leq \frac{1}{i} C \, . 
\]
We would like to use the bundle $F$ to construct a bundle on $M$ with connection 
whose curvature norm is related to $C$ in the same way. We think 
of the fibre over $x \in M$ as the direct sum of all the fibres $F_y$ where $y \in f^{-1}(x)$. 
If the degree of $\overline{M} \to M$ is infinite then this 
``wrapped up'' bundle will have infinite dimension and   
for arriving at the statement of Proposition \ref{technisch} we must keep control 
of the ``structure group'' of this infinite dimensional bundle.

Because $\psi$ is constant outside $K$, there is a unitary trivialization 
\[
          F|_{\overline{M} \setminus K} \cong  (\overline{M} \setminus K)\times \C^d \, .  
\]
such that with respect to this trivialization, the connection $1$-form of the connection 
on $F$ vanishes. 

Moreover, there is a finite cover 
\[
    U_1, \ldots, U_k 
\]
of $M$ consisting of open connected subsets such that 
the restriction of $F$ to $f^{-1}(U_{\kappa})$ is trivial for all $\kappa \in \{1,2, \ldots, k\}$. 
By choosing the members of this cover sufficiently small, we can furthermore assume that 
for all $\kappa ,\lambda \in \{1, \ldots, k\}$ and for each component $V \subset f^{-1}(U_{\kappa})$ 
there is at most one component $W \subset f^{-1}(U_{\lambda})$ with $V \cap W \neq \emptyset$
and that  each $U_{\kappa}$ is evenly covered by the universal cover 
\[ 
   \pi :   \widetilde{M} \to M 
\]
i.e. each component of $\pi^{-1}(U_{\kappa})$ is 
mapped diffeomorphically to $U_{\kappa}$. 

Pick a basepoint in $\overline{M}$ and the induced base point in $M$.  Now the fundamental 
group 
\[
    H := \pi_1(\overline{M}) 
\]
acts freely and transitively from the right 
on $\widetilde{M}$. For each $\kappa$ this induces a parameterization of the set 
of components of $f^{-1}(U_{\kappa})$ by $G/H = \{gH ~|~ g \in G \}$, where 
\[
    G := \pi_1(M)  \, . 
\]
However this parameterization is not canonical, but depends on singling out a
preferred component of $f^{-1}(U_{\kappa})$. We assume that this has been 
done.  
For each $\kappa$ and each component $V \subset f^{-1}(U_{\kappa})$, we choose a unitary 
trivialization of $F|_V$. If $V \cap K = \emptyset$, we take the trivialization 
which was already chosen before. 

Let  
\[   
    \mathcal{H} := l^2(G/H) \otimes \C^d  
\]
be the Hilbert space of square summable $G/H$-families of vectors in $\C^d$. This space
is finite dimensional, if $\overline{M}$ is compact, and a separable infinite dimensional 
Hilbert space, if $\overline{M}$ is not compact ($G/H$ is then countably infinite). 
For each $\kappa, \lambda \in \{1, 2, \ldots, k\}$, we obtain a transition function 
\[
    \alpha_{\kappa \lambda} : (U_{\kappa} \cap U_{\lambda}) \to \mathcal{B}(\mathcal{H})  
\]
describing the passage from the trivialization of $F$ over $f^{-1}(U_{\kappa})$ to 
the trivialization over $f^{-1}(U_{\lambda})$. Because $F$ is unitary, 
the image of $\alpha_{\kappa \lambda}$ is contained in the group of isometries of 
$\mathcal{H}$.

The group of permutations of $G/H$ acts isometrically on $\mathcal{H}$ by precomposition
in $l^2(G/H)$. Let 
\[ 
    C_S \subset \mathcal{B}(\mathcal{H}) 
\]
be the norm closure of the algebra generated by these isometries. Furthermore, 
let 
\[   
    C_T \subset \mathcal{B}(\mathcal{H}) 
\]
be the $C^*$-algebra generated by ``elementary'' linear maps 
\[
        \gamma_0 \otimes \C^d \to \mu_0  \otimes \C^d
        \]
extended by zero on the orthogonal complement of $\gamma_0\otimes \C^d\subset
l^2(G/H)\otimes \C^d$,
where $\gamma_0,\mu_0 \in G/H$ are arbitrary (considered as elements in
$l^2(G/H)$). 
Let 
\[
    C_{S,T} \subset \mathcal{B}( \mathcal{H}) 
\]
be generated by $C_S$ and $C_T$. Note that 
$C_{T} = C_{S,T}$, if $\overline{M}$ is compact, and that  
$C_T$ is a proper subalgebra of $C_{S,T}$ which can be identified with $\K$, if $\overline{M}$ 
is not compact. Furthermore, $C_T$ is a two sided ideal in $C_{T,S}$. 

We now set
\[
    C_i := \{\ (c_1, c_2) \in C_{S,T} \times C_{S,T} ~|~ c_1 - c_2 \in C_{T} \}     \, . 
\]
An easy calculation shows that this is indeed a $C^*$-algebra. It is equipped with a unit. 

\begin{rem} The construction of this algebra is inspired by the
  discussion in John Roe's article \cite{Roe}.
\end{rem}

The projection $C_{S,T} \times C_{S,T} \to C_{S,T}$ onto the second factor 
induces a split exact sequence 
\[
    0 \to C_{T}  \to C_i \to C_{S,T} \to 0 \, , 
\]
the splitting being induced by the diagonal embedding 
\[
   C_{S,T} \to C_{S,T} \times C_{S,T} \, . 
\]
Because the cover $(U_{\kappa})$ is sufficiently small (in the sense explained before), 
for each transition function $\alpha_{\kappa \lambda}$ there is 
a uniquely determined permutation $\beta_{\kappa \lambda}$ of $G/H$ such that the image of the composition 
\[
     (\beta_{\kappa \lambda})^{-1} \circ \alpha_{\kappa \lambda} 
\]
is contained in a $G/H$-family of unitary maps $\C^d \to \C^d$ all but finitely many
of which are equal to the identity. This gives rise to relative transition functions  
\begin{eqnarray*}
    \rho_{\kappa \lambda} : U_{\kappa} \cap U_{\lambda} & \to & C_{S,T} \times C_{S,T} \, , \\
                                        x & \mapsto & (\alpha_{\kappa, \lambda}(x), 
                                                       \beta_{\kappa, \lambda}(x))  \, . 
\end{eqnarray*} 

Note that $\beta_{\kappa,\lambda}$ is locally constant, and that these
transition functions satisfy the cocycle condition. The following lemma is
immediate.  

\begin{lem} \label{smooth} Each relative transition function 
$\rho_{\kappa \lambda}$ is smooth and has image contained in the unitaries of $C_i$. 
\end{lem}

We now set 
\[
    F_i := \prod_{j = 1, \ldots, k} (j, U_{j} \times C_i) / \big(
      (\kappa , x, v) \sim (\lambda, x, \rho_{\kappa \lambda}(x) \cdot v) \big) \, . 
\]
By Lemma \ref{smooth} this is a smooth Hilbert $C_i$-module bundle over $M$. Its fibres are isomorphic 
to $C_i$ and are equipped with the $C_i$-valued inner product induced by
\[
      C_i \times  C_i \to C_i \, , (x,y) \mapsto x^* \cdot y \, . 
\]
Over each subset $U_{\kappa} \subset M$, a section of $F_i|_{U_{\kappa}}$ 
defines a $G/H$-family of sections of the trivial bundle $\underline{\C^d} \to U_{\kappa}$ 
(corresponding to the components of $f^{-1}(U_{\kappa})$). Furthermore, 
the connection on $F$ induces a $G/H$-family of connections on $\underline{\C^d}$ all 
but finitely many of which are trivial. The last property is implied by the choice of trivialization of 
$F|_{\overline{M} \setminus K}$.  These considerations show that we get an induced $C_i$-linear
metric connection 
\[
    \nabla_i : \Gamma(F_i)  \to \Gamma(T^*M \otimes F_i) \, .   
\]
In the above trivialization of $F_i$, this is given by $C_i \cap (C_{T} \times 0)$-valued $1$-forms on each $U_j$. 
Note that over each $U_{\kappa} \subset M$, the curvature 
\[
    \Omega^2(U_{\kappa}; \End_{C_i} (F_i) ) = \Omega^2 (U_{\kappa}; C_i) 
\]
can be viewed as a form in $\Omega^2(U_{\kappa})$ with values 
in a $G/H$-family of matrices $\mathfrak{u}(d)$, all but finitely many of
which are equal to zero. Note also  that the norm 
$\| \Omega_{F_i} \|$ of the curvature of $\nabla_i$ is equal to $\| \Omega_F\|$. In particular, 
it satisfies the inequality
\[
     \| \Omega_{F_i} \| \leq \frac{1}{i} C  
\]
which we have been aiming at. 

If $\overline{M}$ is not compact, then $C_{T}$ can be identified with $\K$ and 
the algebra $C_i$ has the form described in Proposition \ref{technisch}.
However, if $\overline{M}$ is compact, then $C_{T}$ is a finite dimensional matrix 
algebra.  In this case we replace the Hilbert space $\mathcal{H}$ 
by its stabilization $\mathcal{H} \otimes l^2(\N)$ and correspondingly replace 
$C_T$, $C_{S}$, $C_{S,T}$ and $C_i$ by $C_T \otimes \K$, etc. The 
inclusion $\C \hookrightarrow l^2(\N)$ onto the first factor induces canonical 
inclusions $C_T \hookrightarrow C_T \otimes \K$, etc. In particular, if 
\[
    p : l^2(\N) \to \C \subset l^2(\N)
\]
denotes the projection onto the first factor, then 
\[
   q_i := 1 \otimes p \in C_i \otimes \K 
\]
is a projection and the transition functions 
$\rho_{\kappa \lambda}$ have images that are contained in the unitaries of $q_i (C_i \otimes \K)$. 
Replacing the typical fibre $C_i$ by $q_i (C_i \otimes \K)$ in the 
definition of the bundle $F_i$, we then arrive at the statement of 
Proposition \ref{technisch}, if we set $q_i := 1$ in the case of non-compact $\overline{M}$. 
>From now on we work with these modified objects (without change of notation), if $\overline{M}$ is compact. 

\begin{rem} This stabilization process is done only for convenience,
  because we want to make sure that $C_i$ 
splits off a summand $\K$ independently of $\overline{M}$ being
compact or not. This makes sure that  
$K_0( \prod C_i)$ splits off a summand $\prod \Z$.  
\end{rem}

It remains to calculate the component in the $\Z$-summand of the index of the spin Dirac operator on $M$ twisted 
with the Hilbert $C_{i}$-module bundle $F_i$ with connection.

We use the Chern-Weil calculus developed in \cite{Sch}. Let 
\[
   \mathcal{T} \subset C_{T} = \K   
\]
be the subalgebra of trace class operators and set 
\[
 D  : = \{\ (c_1, c_2) \in C_{S,T} \times C_{S,T} ~|~ c_1 - c_2 \in 
          \mathcal{T} \} \subset C_i \,  .
\]
Note that the projection $q_i$ is contained in $D$. We get again a split exact sequence 
\[
    0 \to \mathcal{T}  \to D \to C_{S,T} \to 0 \, . 
\]
The algebra $\mathcal{T}$ is a two sided ideal in $D$.
The sum of the trace-norm on $\mathcal{T}$ and the $C^*$-algebra norm on 
$C_{S,T}$ induce a Banach algebra structure on $D$ so that the inclusion 
\[
     D \to C_i 
\]
is continuous. These facts follow from standard properties of the trace-norm (see
\cite[ I.8.5.6]{Black}) and the triangle inequality for the $C^*$-algebra 
norm on $C_i$. Consider the map 
\begin{eqnarray*}
     \tau :  D & \to & \C \, \\
            (c_1, c_2) & \to & \tr(c_1- c_2) \, , 
\end{eqnarray*}
using the canonical trace on $\mathcal{T}$. It is easily verified that $\tau$ defines 
a trace on $D$ and that 
\[
   \tau(c) = \tr(c_{\mathcal{T}}) \, , 
\]
where $c_{\mathcal{T}}$ is the projection of $c$ onto the summand $\mathcal{T}$ in $D = \mathcal{T} \oplus C_{S,T}$.

In order to apply the results from \cite{Sch}, we need the following lemma. 

\begin{lem} \label{dense} $D$ is a dense subalgebra of $C_i$  and closed under holomorphic 
functional calculus in $C_i$. 
\end{lem} 

\begin{proof} The first assertion is obvious. For the second assertion we need to show 
that if $c$ is an element in the matrix 
algebra $M_k(D)$ and $f$ is a holomorphic function defined in a neighbourhood of 
the spectrum of $c$ regarded as an element in $M_k(C_i)$, then $f(c) \in M_k(D)$. 
By the construction of holomorphic functional calculus via the Cauchy integral 
formula (recall that $D_i$ is a Banach algebra), it is enough to show that each 
element in $D$ which is invertible in $C_i$ is also invertible in $D$
(and the same for matrices). 

Let 
\[
    x = x_1 + x_2 \in  \mathcal{T} \oplus C_{S,T} = D 
\]
be an element which is invertible in $C_i$. Then, by 
assumption, there is an element 
\[
   y = y_1 + y_2  \in \K  \oplus C_{S,T} = C_i 
\]
with 
\begin{eqnarray*}
    (x_1 + x_2) (y_1 + y_2) & = & 1     \, , \\
    (y_1 + y_2) (x_1 + x_2) & = & 1      \, . 
\end{eqnarray*}
This implies that $x_2 y_2 = 1$ and $y_2 x_2 = 1$, 
because $\K$ and $\mathcal{T}$ are $2$-sided ideals in $C_i$. 
Hence, $x_2$ is invertible in $C_{S,T}$ and 
\begin{eqnarray*}
    (x_2^{-1} x_1  + 1) (x_2  y_1 + 1) & = & 1 \, , \\
    (x_2 y_1 +1)(x_2^{-1} x_1 +1) & = & 1 \, .   
\end{eqnarray*}
Regarding $\xi := x_2^{-1} x_1 +1$ as an element in $\mathcal{T}_+$, 
the unitalization of $\mathcal{T}$, this means 
that $\xi$ is invertible in $\K_+$. But because each element 
in $\mathcal{T}_+$ that is invertible in $\mathcal{K}_+$ is already
invertible in $\mathcal{T}_+$, the element $\xi$ is invertible in 
$\mathcal{T}_+$ and therefore $x_2 y_1 +1 \in \mathcal{T}_+$.
We conclude $y_1 + y_2 \in D$.

The same argument applies to matrices over $D$.
\end{proof}

This lemma implies that the inclusion 
\[
    D \to C_i 
\]
induces an isomorphism after applying the $K$-theory functor (see \cite{B}, III.5). 
The trace $\tau:D \to \C$ 
induces a group homomorphism $K_0(D) \to \Z$ that we also denote by $\tau$. Because
the $K$-theories of $D$ and $C_i$ coincide, we may consider $\tau$ as being 
defined on $K_0(C_i)$. The restriction of this homomorphism to 
the canonical  $\Z$-summand in $K_0(C_i)$ is the identity. 

Each of the transition functions $\rho_{\lambda \kappa}$ takes values 
in the unitaries of $q_i D$. Hence, similarly to the definition of $F_i$, 
we get a smooth Hilbert $D$-module bundle 
$\mathcal{W} \to M$ with connection $\nabla_{\mathcal{W}}$. 

We are now in a situation where Theorem 9.2 of \cite{Sch} can be applied. 
In our case it says that  (up to sign) 
\[
   z_i =  \tau(\ind(D_{F_i})) =  \int_M  \mathcal{A}(M) \cup \ch_{\tau}(\mathcal{W},\nabla_{\mathcal{W}}) \, . 
\]
For the definition of $\ch_{\tau}$, we refer to \cite{Sch}, Definition 5.1.

\begin{rem} The number $z_i$ can be determined without Lemma \ref{dense} 
by calculating first $\tau(\ind(D_{\mathcal{W}}))$ and applying a diagram chase 
comparing the split exact $K$-theory sequences for the algebras $D$ and $C_i$. 
However, the first step requires an appropriate (yet unproblematic) 
reformulation of the Chern-Weil calculus of \cite{Sch}. 
The above argument avoids this and fits more neatly the statement 
of Theorem 9.2 in \cite{Sch}. 
\end{rem}

Let us calculate the right hand side over an open subset $U_{\kappa} \subset M$. Using the 
trivialization of $F|_{f^{-1}(U_{\kappa})}$ from above and  the definition of
the trace $\tau$, we have 
\[
    \int_{U_{\kappa}} \mathcal{A}(M) \cup \ch_{\tau}(\mathcal{W},\nabla_{\mathcal{W}}) = \int_{f^{-1}(U_{\kappa})} f^*(\mathcal{A}(M)) 
               \cup (\ch(F - \underline{\C^{d}},\nabla_{F}-\nabla_{\C^d})) \, . 
\]
Note that the right hand integral involves only finitely many components of $f^{-1}(U_{\kappa})$. 
Hence, by use of a partition of unity subordinate to the cover $(U_{\kappa})$ of $M$, we get
\[
   z_i = \tau(\ind(D_{F_i})) = \int_{\overline{M}} \mathcal{A}(\overline{M}) \cup (\ch(F - \underline{\C^d},\nabla_F-\nabla_{\C^d})) \, .
\]
By construction, $\ch(F - \underline{\C^d})$ represents a compactly supported nonzero cohomology class 
of degree $n$ in $\overline{M}$ . This implies $z_i \neq 0$ and completes  
the proof of Proposition \ref{technisch}. 

Note that, by the relative index formula of \cite{GL}, $z_i$ is exactly the
relative index of Dirac operator on $\overline M$ twisted on the one hand with
the trivial bundle of rank $d$, and on the other hand twisted with $F$. In
retrospective, we gave a definition of the relative index using suitable
$C^*$-algebras and their K-theory (as in \cite{Roe}) and reduced the proof of
the corresponding index formula to the methods described in \cite{Sch}.

\end{document}